\documentclass[titlepage,11pt]{article}

\usepackage[margin=1.2in]{geometry}
\usepackage{latexsym}
\usepackage{amsfonts,amsmath}

\usepackage{tikz}
\usepackage{float}
\usepackage{url}
\usetikzlibrary{arrows}
\usetikzlibrary{decorations.markings}
\usepackage[symbol]{footmisc}

\usepackage{xspace}
\usepackage[color=Olive]{todonotes}

\usepackage[colorlinks=true,allcolors=blue]{hyperref}
\usepackage{tikz}
\usetikzlibrary{decorations.pathmorphing}
\tikzset{every node/.style={circle,fill, inner sep=0pt, minimum size=1.5mm}, every picture/.style={scale=.85}}

\usepackage{hyperref, pythonhighlight}

\usepackage{lineno}
\modulolinenumbers[5]

\usepackage{amsmath,amssymb,amsthm}
\allowdisplaybreaks
\renewcommand{\epsilon}{\varepsilon}
\newcounter{tbox}

\newtheorem{lemma}{Lemma}

\newtheorem{conjecture}{Conjecture}

\newtheorem{claim}{Claim}

\newtheorem{thm}{Theorem}

\usepackage{thm-restate}

\usepackage{authblk, url}

\title{Improved bounds for the triangle case of Aharoni's rainbow generalization of the Caccetta-H\"{a}ggkvist conjecture\thanks{This is an accepted manuscript. The published version appeared in Discrete Mathematics, 
Volume 347, Issue 1, January 2024, 113691, \url{https://doi.org/10.1016/j.disc.2023.113691}}}

\author[1]{Patrick Hompe}
\author[1]{Zishen Qu\thanks{Current address: Department of Mathematics
University of Illinois at Urbana-Champaign, Urbana, IL.}}
\author[1]{Sophie Spirkl\thanks{We acknowledge the support of the Natural Sciences and Engineering Research Council of Canada (NSERC), [funding reference number RGPIN-2020-03912]. Cette recherche a été financée par le Conseil de recherches en sciences naturelles et en génie du Canada (CRSNG), [numéro de référence RGPIN-2020-03912].  This project was funded in part by the Government of Ontario.
}}

\affil[1]{University of Waterloo, Department of Combinatorics and Optimization}

\begin{document}
\maketitle
		
		\begin{abstract}
For a digraph $G$ and $v \in V(G)$, let $\delta^+(v)$ be the number of out-neighbors of $v$ in $G$. The Caccetta-H\"{a}ggkvist conjecture states that for all $k \ge 1$, if $G$ is a digraph with $n = |V(G)|$ such that $\delta^+(v) \ge k$ for all $v \in V(G)$, then $G$ contains a directed cycle of length at most $\lceil n/k \rceil$. Aharoni proposed a generalization of this conjecture, that a simple edge-colored graph on $n$ vertices with $n$ color classes, each of size at least $k$, has a rainbow cycle of length at most $\lceil n/k \rceil$. Let us call $(\alpha, \beta)$ \emph{triangular} if every simple edge-colored graph on $n$ vertices with at least $\alpha n$ color classes, each with at least $\beta n$ edges, has a rainbow triangle. Aharoni, Holzman, and DeVos showed the following:

\begin{itemize}

  \item $(9/8,1/3)$ is triangular;
  \item $(1,2/5)$ is triangular.

\end{itemize}
In this paper, we improve those bounds, showing the following:
 
 \begin{itemize}

  \item $(1.1077,1/3)$ is triangular;
  \item $(1,0.3988)$ is triangular.
\end{itemize}
Our methods give results for infinitely many pairs $(\alpha, \beta)$, including $\beta < 1/3$; we show that $(1.3481,1/4)$ is triangular.
\end{abstract}

\section{Introduction and preliminaries}
We call a graph $G$ \emph{simple} if it has no loops or parallel edges, and let a digraph $D$ is a \emph{simple} digraph if the underlying undirected graph has no loops or parallel edges. For a digraph $D$ and a vertex $v \in V(D)$, we let $\delta^+(v)$ denote the number of out-neighbors of $v$ in $D$. The following conjecture was made by Caccetta and H\"{a}ggkvist in 1978, and despite much study still remains open:

\begin{conjecture}[\cite{ch}]\label{ch_thm}
Let $n,k$ be positive integers. Suppose that $G$ is a simple digraph on $n$ vertices with $\delta^+(v) \ge k$ for all $v \in V(G)$; then $G$ contains a directed cycle of length at most $\lceil n/k \rceil$.
\end{conjecture}

For a graph $G$, an \emph{edge-coloring} of $G$ is an assignment of a color to each edge in $E(G)$, which we denote by $f: E(G) \rightarrow \mathbb{N}$. Note that we do not require this to be a proper edge-coloring. For a graph $G$ with an edge-coloring $f$, a \emph{rainbow cycle} in $G$ (with respect to $f$) is a cycle $C$ in $G$ such that for all pairs of distinct edges $e_1,e_2 \in E(C)$, we have that $f(e_1) \ne f(e_2)$. For $i \in \mathbb{N}$, the set of edges $f^{-1}(i)$ is called a \emph{color class} (corresponding to color $i$). The \emph{rainbow girth} of an edge-colored graph $G$ is the length of a shortest rainbow cycle in $G$, and $\infty$ if not such cycle exists. 

In \cite{aharoni}, Aharoni proposed a generalization of Conjecture \ref{ch_thm}:

\begin{conjecture}[\cite{aharoni}]\label{aharoni_thm}
Let $n,k$ be positive integers. Suppose that $G$ is a simple graph on $n$ vertices. Let $b$ be a coloring of the edges of $G$ with $n$ color classes of size at least $k$; then $G$ has rainbow girth at most $\lceil n/k \rceil$.
\end{conjecture}

Recently, there have been a number of results related to this conjecture. In \cite{devos}, the conjecture is proved for color classes of size equal to $2$, which is generalized in \cite{chudnov} to a statement for the case where the color classes have size equal to one or two. In \cite{hompe}, Conjecture \ref{aharoni_thm} is proven for color classes of size at least $ck$ for $c = 10^{11}$. In \cite{guo}, the case where the color classes each form a matching is investigated. Finally, in \cite{pt}, Conjecture \ref{aharoni_thm} is proved up to an additive constant; that is, it is shown that for every $k$, there is a constant $c_k$ such that under the assumptions of Conjecture \ref{aharoni_thm}, $G$ has rainbow girth at most $n/k + c_k$. 

For the Caccetta-H\"{a}ggkvist conjecture, there is a series of papers investigating the ``triangle case,'' namely that a simple digraph on $n$ vertices with minimum out-degree at least $n/3$ contains a directed triangle. This seems to be the hardest case and perhaps characteristic of the problem as a whole. Also, the existence of rainbow triangles in graphs has been well-studied in the literature (under different conditions, see \cite{aharoni}). In this paper, we focus on the triangle case of Aharoni's rainbow conjecture, which states that if there are $n$ color classes of size at least $n/3$, then there is a rainbow triangle. In the triangle case, a cycle is rainbow if and only if it is properly colored, meaning that no two adjacent edges in the cycle have the same color. This makes, in the triangle case, the rainbow conjecture equivalent to a similar conjecture made in \cite{aharoni} about the existence of short properly-colored cycles.  Let us say that $(\alpha, \beta)$ is \emph{triangular} if every simple edge-colored graph on $n$ vertices with at least $\alpha n$ color classes, each with at least $\beta n$ edges, has a rainbow triangle. Now, in \cite{aharoni}, the following partial results in this direction are shown:
\begin{thm}[\cite{aharoni}]\label{existing}
$(9/8,1/3)$ is triangular; $(1,2/5)$ is triangular.
\end{thm}

We improve both of the above results, and we show a family of results for the case where the color classes have size less than $n/3$. We present an interesting example from this family along with the two improvements in the following two theorems:

\begin{restatable}{thm}{mainone}\label{main1}
$(1.1077,1/3)$ is triangular; $(1.3481,1/4)$ is triangular.
\end{restatable}

\begin{restatable}{thm}{maintwo}\label{main2}
	$(1,0.3988)$ is triangular.
\end{restatable}

An outline of the paper is as follows. In Section 2, we prove some auxiliary lemmas. In Section 3, we consider the case of $(1+\delta)n$ colors for $\delta > 0$, showing Theorem \ref{main1}, and in fact a more general result. Then, in Section 4, we consider the case of $n$ color classes each of size at least $\beta n$, and show Theorem \ref{main2}. In the final section of the paper, we make some concluding remarks, and suggest areas of further work. In our proof, we make use of the following result (which we will also modify to obtain a stronger bound under certain conditions), due to Goodman, which was used in \cite{aharoni} for the proof of Theorem \ref{existing} as well.

\begin{thm}[\cite{goodman}]\label{goodman_thm}
Suppose $G$ is a simple graph with $n$ vertices and $m$ edges, and let $t(G)$ denote the number of triangles in $G$; then we have:
$$t(G) \ge \frac{4m}{3n}\left(m-\frac{n^2}{4}\right) $$
\end{thm}

We also use the following state-of-the-art approximate result on the Caccetta-H\"{a}ggkvist conjecture:

\begin{thm}[\cite{ch_result}]\label{ch_result_thm}
Let $D$ be a simple digraph on $n$ vertices with minimum out-degree at least $0.3465n$; then $D$ contains a directed triangle.
\end{thm}

We note that improvements to the approximation in Theorem \ref{ch_result_thm} would result in improvements to Theorem \ref{main2}.

\section{Counting triples}

Let us first prove the following straight-forward fact: 
\begin{lemma} \label{lem:convex}
    Let $k \in \mathbb{R}$, and let $$f_k(x) = {x \choose 2} + {k-x + 1 \choose 2} = x(x-1)/2 + (k-x+1)(k-x)/2.$$
    Then $f_k(x)$ is convex with a unique global minimum at $x = (k+1)/2$.
\end{lemma}
\begin{proof}
    We observe that $f'_k(x) = x - \frac{1}{2} -(k-x) - \frac{1}{2}$ vanishes when $2x = k+1$ so $x = (k+1)/2$. Moreover, $f_k''(x) = 2$, and so $f_k$ is convex with a unique global minimum at $x = (k+1)/2$. 
\end{proof}

Let $G$ be a graph. We call a set $\{u, v, w\} \subseteq V(G)$ a \emph{happy triple} if $G[\{u, v, w\}]$ has at least two edges. We consider the maximum number of happy triples in graphs with bounded degree. A similar question (maximizing the number of two-edge paths in graphs with a bounded number of vertices) was considered in \cite{ak}. 
\begin{lemma} \label{lem:ht} Let $k, l \in \mathbb{N}$. Let $G$ be a graph with $k$ edges and with maximum degree at most $l \geq k/2$. Then $G$ contains at most ${l \choose 2} + {k-l+1 \choose 2}$ happy triples. 
\end{lemma}
\begin{proof}  
    Fix $k$ and $l \ge k/2$ with $k,l \in \mathbb{N}$. For every two distinct edges $e, f$ in $E(G)$, either $e \cup f$ is a happy triple, or $|e \cup f| = 4$, and so $G$ has at most ${k \choose 2}$ happy triples. Therefore, we may assume that $l < k$. 
    
    Let $G$ be a graph that satisfies the assumptions of the lemma, and subject to that, has as many happy triples as possible. We may assume that $G$ has more than ${l \choose 2} + {k-l+1 \choose 2}$ happy triples. By Lemma \ref{lem:convex}, this quantity is minimized when $l = (k+1)/2$; so we may assume that $G$ has more than $2{k/2+1/2 \choose 2} = (k/2-1/2)(k/2+1/2) = k^2/4 - 1/4$ happy triples. Because each happy triple involves at least two edges, this implies that on average, each edge of $G$ is part of more than $(k-1/k)/2$ happy triples, and so one of its ends has degree more than $(k-1/k)/4 + 1 > k/4$. 

    Let us call a vertex \emph{important} if it has degree more than $k/4$. The above shows that $G$ contains an important vertex. By counting edges, we get that the number $s$ of important vertices satisfies $sk/4 < 2k$ and so $s < 8$. 
    
     \begin{claim}\label{claim:degmax}
        The graph $G$ does not contain a vertex of degree at least $k/2$. 
    \end{claim}

    Let $v$ be a vertex of degree $l' \geq k/2$ in $G$. We have ${l' \choose 2}$ happy triples that contain two edges incident with $v$. For every edge $xy$ of $G$ not incident with $v$, there is at most one happy triple $\{x, y, v\}$, and so there are at most $k-l'$ happy triples that involve $v$ but do not contain two edges incident with $v$. 

    Moreover, the number of happy triples not involving $v$ is the number of happy triples formed by choosing two edges not incident with $v$, which is at most ${k - l' \choose 2}$. This leads to the overall bound of ${l' \choose 2} + {k-l' \choose 2} + k-l' = {l' \choose 2} + {k - l' + 1 \choose 2}$ on the number of happy triples. It remains to show that this is maximized when $l = l'$.  
    By Lemma \ref{lem:convex}, the function $f_k$ is convex, and so $f_k$ is maximized at one of the ends of the range: either at $l' = l$ (and we are done) or at $l' = \lceil k/2 \rceil$. If $k$ is odd, then $\lceil k/2 \rceil = (k+1)/2$ is the global minimum of $f$ by Lemma \ref{lem:convex}. If $k$ is even, we notice that $f_k(k/2) = f_k(k/2+1)$; so either $l \in \{k/2, k/2+1\}$ and thus $f_k(l) = f_k(l')$, or $l > k/2+1$ and $f_k(l) > f_k(l')$ since $f'_k(x) > 0$ for all $x > k/2 + 1$. This proves Claim \ref{claim:degmax}. 
    
\medskip
    It remains to consider the case that all vertices of $G$ have degree less than $k/2$. In particular, this means that if we can move one edge of $G$ and increase the number of happy triples, then we have found a contradiction, because moving one edge increases the maximum degree by at most 1, and so the maximum degree is still at most $l$ after moving one edge. 
    
    Let $S$ be the set of vertices of $G$ of degree more than $k/8$. Then, $|S|k/8 < 2k$ by counting edges, and so $|S| < 16$. 
    
    Suppose that $ab \in E(G)$ with $a, b \not\in S$. Then $ab$ is part of at most $k/4$ happy triples in $G$. Let $v$ be an important vertex of $G$. Then, if we consider the graph $G'$ arising from $G$ by removing $ab$, and adding a vertex $z$ and an edge $vz$, it follows that $G'$ has more happy triples than $G$, because $vz$ is now part of more than $k/4$ happy triples. 

    It follows that $S$ is a vertex cover of $G$. Let $S = \{s_1, \dots, s_t\}$ where $\deg(s_1) \geq \dots \geq \deg(s_t)$ and $t < 16$. We may assume that $G$ is chosen with $\deg(s_1)$ maximized (subject to satisfying the constraints of the lemma, and maximizing the number of happy triples), and if $\deg(s_1) \geq k/2$, then the claim holds, so we may assume that $\deg(s_1) < k/2$. Also, recall that $G$ contains an important vertex, that is, a vertex of degree more than $k/4$. Since $s_1$ is a vertex of maximum degree in $G$, it follows that $\deg(s_1) > k/4$. 
    
    Let $r \in \mathbb{N}$ be minimal such that $S' = \{s_1, \dots, s_{r}\}$ is a vertex cover of $G$. Then $r > 1$, because otherwise $s_1$ has degree $k$; but we assumed that $\deg(s_1) < k/2$. 
    
    From the choice of $r$, it follows that $G$ contains an edge $s_rx$ for some $x \not\in S'$. Then, $s_rx$ is part of at most $\deg(s_r)-1 + r-1$ happy triples, because $\deg(x) \leq r$ (using that $x \not\in S'$, and $S'$ is a vertex cover, so $N(x) \subseteq S'$). Let $G'$ be the graph arising from $G$ by removing $s_rx$ and adding a new vertex $z$ and edge $s_1z$. This creates $\deg(s_1)$ new happy triples, and so from the choice of $G$ with $\deg(s_1)$ maximized, we may assume that $G'$ has strictly fewer happy triples than $G$ does. Therefore, $\deg(s_1) < \deg(s_r) - 1 + r -1 = \deg(s_r) + r - 2$, and since degrees are integers, it follows that $\deg(s_1) \leq \deg(s_r) + r - 3$. 

    This immediately implies that $r \geq 3$. By counting edges, we have 
    \begin{align*}
        k + {r \choose 2} &\geq \sum_{i=1}^r \deg(s_i)  \geq r \deg(s_1) - (r-3)(r-1). 
    \end{align*}

    So $\deg(s_1) \leq \frac{k}{r} + \frac{3}{2}r - \frac{7}{2}$. Using that $s_1$ is important, that is, $\deg(s_1) > k/4$, this implies that $k/r + 3r/2 - 7/2 > k/4$, so $3r/2 - 7/2 > (r-4)k/(4r)$. Therefore, for $r \geq 5$, this gives $k < \frac{4r(3r/2 - 7/2)}{r-4}$.   This implies that $k \leq 103$ (by plugging in possible values of $r \in \{5, \dots, 15\}$) when $r \geq 5$. 

    We now address the case $r \le 4$. For each happy triple $\{x, y, z\}$, we choose a vertex in $\{x, y, z\}$ which has two neighbors in $\{x, y, z\}$ and call this the \emph{center} of $\{x, y, z\}$. 

    The number of happy triples is at most 
    $$\frac{r}{2}\left(\frac{k}{r} + \frac{3}{2}r - \frac{7}{2}\right)\left(\frac{k}{r} + \frac{3}{2}r - \frac{9}{2}\right) + k(r-1)/2.$$
    The former term is an upper bound for happy triples with center in $S'$, and the latter is an upper bound for happy triples with center outside $S'$: each such triple consists of an edge $s_ix$ for some $s_i \in S'$, and a vertex $s_j \in S' \setminus \{s_i\}$; and counting this way counts each such triple twice (once for each way to choose $i$ and $j$). 

    So we have 
    \begin{align} \label{eq:k}
        \frac{r}{2}\left(\frac{k}{r} + \frac{3}{2}r - \frac{7}{2}\right)\left(\frac{k}{r} + \frac{3}{2}r - \frac{9}{2}\right) + k(r-1)/2 &> k^2/4 - 1/4.
    \end{align}
    For each fixed value of $r$, this is a quadratic in $k$, and by plugging in the possible values of $r \in \{3, 4\}$, we get that $k \leq 30$ (explicit calculations are shown in \ref{ap:1}), and so, combined with our above-mentioned bound for $r \geq 5$, we obtain that $k \leq 103.$ This leaves a finite number of cases, and we handle them using dynamic programming, see \ref{ap:2}. 
\end{proof}

We briefly remark that the bound in Lemma \ref{lem:ht} is best possible: When $l > k/2$, it is achieved by taking $K_{1,l-1}$ and $K_{1, k-l}$ and adding an edge between their centers; and then for $l = k/2$, it is achieved by $K_{2, l}$. It would be interesting to consider the best upper bound for the regime when $l < k/2$. Our program does not give much insight into this; for example, for $l = 3$ and $k = 3t$, the value output by our program is for $K_{t, 3}$, which does not actually have maximum degree 3 for $t \ge 4$. 

Next, we refine Theorem \ref{goodman_thm}, as follows:
\begin{thm}\label{triange_count}
    Let $H$ denote the unique graph with three vertices and one edge. Let $G$ be a simple graph on $n$ vertices with $m$ edges, and let $h(G)$ denote the total number of induced copies of $H$ in $G$. Then the number of triangles in $G$ is at least:
    $$ \frac{h(G)}{3} + \frac{4m}{3n}\left(m-\frac{n^2}{4}\right).$$
\end{thm}
\begin{proof}
For a vertex $v \in V(G)$, we let $N(v) \subseteq V(G)$ be the set of neighbors of $v$ (not including $v$), and we let $\deg(v) = |N(v)|$. Recall that for two sets of vertices $A, B \subseteq V(G)$, we have by basic set properties that $|A \cup B| = |A| + |B| - |A \cap B|$. Now, for each edge $e = uv \in E(G)$, if we let $A_u = N(u) \setminus \{v\}$ and $A_v = N(v) \setminus \{u\}$, we have that the number of triangles containing $e$ is equal to:
$$|A_u \cap A_v| = |A_u| + |A_v| - |A_u \cup A_v|.$$

Let $h_e$ denote the number of copies of $H$ in $G$ containing $e$. Note that $|A_u \cup A_v| \le n - h_e -2$. It follows that the number of triangles which $e$ is contained in is at least $|A_u| + |A_v| - |A_u \cup A_v| \geq |A_u| + |A_v| - n + h_e + 2 = \deg(u) + \deg(v) - n + h_e$.

Now, let $t(G)$ denote the number of triangles in $G$. Summing this over all edges $e \in E(G)$, and noting that we count each triangle at most three times, we obtain:
\begin{align*}
    t(G) &\ge \frac{1}{3} \sum_{e=(uv) \in E(G)} (\deg(u) + \deg(v) - n + h_e) \\
    &= \frac{h(G)}{3} - \frac{mn}{3} + \frac{1}{3}\sum_{v \in V(G)} \deg(v)^2,
\end{align*}
and by the Cauchy-Schwarz inequality, we obtain that this is at least:
\begin{align*}
    \frac{h(G)}{3}-\frac{mn}{3}+\frac{4m^2}{3n} = \frac{h(G)}{3} + \frac{4m}{3n}\left(m-\frac{n^2}{4}\right),
\end{align*}
as desired. This completes the proof.\end{proof}

\section{$(1+\delta)n$ colors}
In this section, we consider the case where we have at least $(1+\delta)n$ color classes each of size at least $tn$. 

We show the following general result:

\begin{thm}\label{main1general}
Let $t,\delta$ be positive real numbers. Suppose $G$ is a simple edge-colored graph with at least $(1+\delta)n$ color classes of size at least $tn$. Then there exists a rainbow triangle if there exists $0 < \epsilon < 1/2$ such that all of the following conditions hold:
\begin{align}
    \frac{(1+\delta-\epsilon \delta)t}{2}<\frac{4}{3}(1+\delta)\left((t(1+\delta)-\frac{1}{4}\right); \label{eq:211} \\
    \frac{8}{3t}\left(t(1+\delta)-\frac{1}{4}\right)(1+\delta)+\frac{t}{12}\left(\lambda+1\right)\left(4(1-\epsilon)-\left(1-\lambda\right)^2\right) >1+\delta
    \label{eq:212};\\
    \frac{16}{3}\left(1+\delta\right)>\frac{2}{3t}+1;\label{eq:213}
\end{align}
where $\lambda = \sqrt{1-2\epsilon}$.
\end{thm}
\begin{proof}
It suffices to show the claim when each color class has size equal to $\lceil tn \rceil$ and when the number of colors is equal to $\lceil (1+\delta)n \rceil$, so we assume both of those conditions. Let $\delta_1, t_1$ be real numbers such that $(1+\delta_1)n = \lceil (1+\delta)n \rceil$ and $t_1 n = \lceil tn \rceil$; so $\delta_1 \geq \delta$ and $t_1 \geq t$. Suppose that there are no rainbow triangles. Then, for each triangle in $G$, at least two of its edges have the same colour. Now, we say that a color class $c$ is \emph{good} if there are at least $(1-\epsilon)\frac{n^2}{2}t_1^2$ triangles in $G$ with at least two of its three edges having color $c$. We first show the following claim.

\begin{claim}\label{first_claim}
More than $n$ of the color classes are good.
\end{claim}
Suppose not; then at most $n$ of the colors are good. We have that at least $\delta_1 n$ of the colors are not good, namely that at least $\delta_1 n$ of the colors have at most $(1-\epsilon)\frac{n^2}{2}t_1^2$ triangles with at least two of the three edges having color $c$. For every triangle $T$ in $G$, at least two of the edges in $T$ have the same color. Letting $t(G)$ denote the number of triangles in $G$, it follows that:
\begin{align*}
    t(G) &\le n\binom{t_1 n}{2} + n\delta_1 (1- \epsilon) \frac{t_1^2 n^2}{2} \\
    &\le n\frac{t_1^2 n^2}{2} + n\delta_1 (1- \epsilon) \frac{t_1^2 n^2}{2} \\
    &= n^3(1+\delta_1 - \epsilon \delta_1)\frac{t_1^2}{2}.
\end{align*}
By Theorem \ref{goodman_thm}, since $m = (1+\delta_1)t_1 n^2$, we have that:
\begin{align*}
    t(G) &\ge \frac{4m}{3n}\left(m-\frac{n^2}{4}\right) = \frac{4(1+\delta_1)t_1 n^2}{3n}\left((1+\delta_1)t_1 n^2-\frac{n^2}{4}\right)\\ &= n^3 \cdot \frac{4(1+\delta_1)t_1}{3}\left((1+\delta_1)t_1-\frac{1}{4}\right).
\end{align*}
Putting together the upper and lower bounds for $t(G)$ and dividing by $n^3t_1$, we conclude that
\begin{align}
    \frac{(1+\delta_1-\epsilon \delta_1)t_1}{2}\geq\frac{4}{3}(1+\delta_1)\left(t_1(1+\delta_1)-\frac{1}{4}\right). \label{eq:5}
\end{align}

Now \eqref{eq:211} asserts that
$$\frac{4}{3}(1+\delta)\left(t(1+\delta)-\frac{1}{4}\right) - \frac{(1+\delta-\epsilon \delta)t}{2} > 0.$$
We claim that for all $t_1 \ge t$ and $\delta_1 \ge \delta$, we have the following analogous inequality, which will lead to a contradiction with \eqref{eq:5}: 
\begin{align}\frac{4}{3}(1+\delta_1)\left(t_1(1+\delta_1)-\frac{1}{4}\right) - \frac{(1+\delta_1-\epsilon \delta_1)t_1}{2} > 0.\label{eq:6} \end{align}

To show this, it suffices to show that the derivative of the above expression is non-negative at all points $(t_1,\delta_1)$ with $t_1 \ge t$ and $\delta_1 \ge \delta$. Indeed, the derivative of the above expression with respect to $t_1$ at the point $(t_1,\delta_1)$ is:
\begin{align*}
    \frac{4(1+\delta_1)^2}{3}-\frac{1+\delta_1-\epsilon \delta_1}{2} \geq \frac{4}{3}(1+\delta_1)^2 - \frac{1}{2}(1+\delta_1) \geq \frac{5}{6}(1+\delta_1)^2> 0
\end{align*}
since $1+\delta_1 \geq 1$. 

Taking the derivative with respect to $\delta_1$ at the point $(t_1,\delta_1)$, we would like to show that: 
\begin{align}
    \frac{8t_1(1+\delta_1)-1}{3}-\frac{(1-\epsilon)t_1}{2} > 0 \label{eq:4}
\end{align}

Now, by \eqref{eq:213}, we have that:
\begin{align*}
    \frac{16}{3}(1+\delta_1) &\ge \frac{16}{3}(1+\delta) 
    > \frac{2}{3t}+1 \\
    &\ge \frac{2}{3t_1}+1 > \frac{2}{3t_1} + (1-\epsilon).
\end{align*}
Now, multiplying by $t_1$ and dividing by 2 implies \eqref{eq:4}. This proves \eqref{eq:6}, and together with \eqref{eq:5}, we obtain a contradiction. This proves Claim \ref{first_claim}. 

\medskip

So, by Claim \ref{first_claim}, we have that more than $n$ colors are good. Next, we show the following.

\begin{claim}\label{alpha}
Suppose that color $c$ is good. Then there exists a vertex in $G$ incident to at least $\alpha n$ edges of $c$, where $\alpha = \frac{t_1}{2}(\sqrt{1-2\epsilon}+1)$.
\end{claim}
Suppose that color $c$ is good and each vertex of $G$ is incident to less than $\alpha n$ edges of $c$. We count the number of triangles in $G$ with at least two edges of color $c$. This is exactly the number of happy triples in a graph $G$ with $t_1n$ edges and maximum degree at most $\lfloor \alpha n \rfloor$. 

We consider two cases. If $\lfloor \alpha n \rfloor \geq t_1n/2$, then let us set $q = \lfloor \alpha n \rfloor$. From the choice of $\alpha$, we have that $\alpha n > t_1n/2$, and so the only remaining case is that $t_1n$ is odd, and $\lfloor \alpha n \rfloor = (t_1n-1)/2$. In this case, let us set $q = (t_1n+1)/2$, which is an integer. In each case, we have that $q \geq t_1n/2$ is an integer and an upper bound for the maximum degree of $G$. 

By Lemma \ref{lem:ht}, the number of happy triples in $G$ is at most 
\begin{align}\label{eq:f}
    f_{t_1n}(q) = \binom{q}{2}+\binom{t_1n-q + 1}{2}
\end{align}
where $f_{t_1n}$ is the function from Lemma \ref{lem:convex}.

Since $c$ is good, we have that the number of happy triples of color $c$ is at least $(1-\epsilon)\frac{n^2 t_1^2}{2}$, so we get that:
\begin{align*}
    (1-\epsilon)\frac{n^2 t_1^2}{2} &\le
    \binom{q}{2}+\binom{t_1n-q + 1}{2} \\
    &= \frac{q^2-q}{2} + \frac{(t_1n-q)^2 + t_1n - q}{2}\\
    &= \frac{q^2}{2} + \frac{(t_1n-q)^2}{2} + \frac{t_1n - 2q}{2}\\
    &\leq \frac{q^2}{2} + \frac{(t_1n-q)^2}{2}.
\end{align*}
using that $q \geq t_1n/2$ in going from the third to the fourth line. 

As a function of $q$, we have that $f(q) = \frac{q^2}{2} + \frac{(t_1n-q)^2}{2}$ is convex and minimized when $q = t_1n/2$. Therefore, if $\lfloor \alpha n \rfloor \geq t_1n/2$, then $q \leq \alpha n$, and so $f(q) \leq f(\alpha n)$, and so, by cancelling a factor $n^2$ from the inequality $$(1-\epsilon)\frac{n^2 t_1^2}{2} \leq f(\alpha n) = \frac{(\alpha n)^2}{2} + \frac{(t_1n-\alpha n)^2}{2},$$ 
we get that 
\begin{align} \label{eq:epseq}
    (1-\epsilon)t_1^2 < \alpha^2 + (t_1-\alpha)^2.
\end{align}

If $q = (t_1n+1)/2$, then $x = q$ is the global minimum of $f_{t_1n}(x)$ by Lemma \ref{lem:convex}. Therefore, the number of happy triples is at most $f_{t_1n}(q) \leq f_{t_1n}(\alpha n)$, so again using that $c$ is good, we have proved that
   \begin{align*}
    (1-\epsilon)\frac{n^2 t_1^2}{2} &\le
    \binom{\alpha n}{2}+\binom{t_1n-\alpha n + 1}{2} \\
    &= \frac{(\alpha n)^2-\alpha n}{2} + \frac{(t_1n-\alpha n)^2 + t_1n - \alpha n}{2}\\
    &= \frac{(\alpha n)^2}{2} + \frac{(t_1n-\alpha n)^2}{2} + \frac{t_1n - 2\alpha n}{2}\\
    &\leq \frac{(\alpha n)^2}{2} + \frac{(t_1n-\alpha n)^2}{2}.
\end{align*}
By cancelling a factor of $n^2$, we again get that \eqref{eq:epseq} holds. 

Now, $\alpha = \frac{t_1}{2}(\sqrt{1-2\epsilon}+1)$ implies that:
\begin{align*}
    \alpha^2 + (t_1-\alpha)^2 &= \left(\frac{t_1}{2}\left(1+\sqrt{1-2\epsilon}\right)\right)^2 + \left(\frac{t_1}{2}\left(1-\sqrt{1-2\epsilon}\right)\right)^2 \\
    &= \frac{t_1^2}{4}\left(2+2(1-2\epsilon)\right) \\
    &= (1-\epsilon)t_1^2 \\
    &< \alpha^2 + (t_1-\alpha)^2,
\end{align*}
using \eqref{eq:epseq}, which is a contradiction, as desired. This proves Claim \ref{alpha}. 

\medskip

Now, since more than $n$ colors are good, Claim \ref{alpha} implies that there exists a vertex $v \in V(G)$ such that there exist two colors $c_1$ and $c_2$ and there are at least $\alpha n$ edges of each of the colors $c_1$ and $c_2$ incident to $v$, where $\alpha = \frac{t_1}{2}(\sqrt{1-2\epsilon}+1)$. 

For $i \in \{1, 2\}$, let $S_i$ be the set of vertices $x$ such that $vx$ has color $c_i$; then $|S_i| \geq \lceil \alpha n \rceil$. 

Now, since $c_1$ is good, there are at least $(1-\epsilon)\frac{n^2 t_1^2}{2}$ triangles in $G$ containing two edges of color $c_1$. The number of such triangles that do not consist of $v$ and two vertices in $S_1$ is at most $\binom{t_1 n - |S_1|}{2} +  (t_1 n - |S_1|)$; the former counts triangles using two edges of color $c_1$ not incident with $v$, while the latter counts triangles of the form $\{v, x, y\}$ with $xy$ of color $c_1$. 

It follows that the number of triangles consisting of $v$ and two vertices in $S_1$, and therefore the number of edges with both ends in $S_1$, is at least:
\begin{align*}
   & \ \ \ \  (1-\epsilon)\frac{n^2 t_1^2}{2} - \binom{t_1 n - |S_1|}{2} - (t_1 n - |S_1|) \\
   &= (1-\epsilon)\frac{n^2 t_1^2}{2}-\frac{(t_1 n  - |S_1|)^2}{2} - \frac{t_1 n - |S_1|}{2} \\
    &\ge (1-\epsilon)\frac{n^2 t_1^2}{2}-\frac{(t_1 n  - \lceil \alpha n \rceil)^2}{2} - \frac{t_1 n - \lceil \alpha n \rceil}{2},
\end{align*}
and since $c_2$ is good, an identical argument gives that the same is true for $S_2$. 

Now, since there is no rainbow triangle in $G$, it follows that edges with one end in $S_1$ and the other in $S_2$ have color $c_1$ or $c_2$. Therefore, there are at most $2 t_1 n - |S_1|-|S_2|$ edges between $S_1$ and $S_2$. 

Let $H$ be the graph on three vertices with exactly one edge. For each edge $e = xy$ with both ends in $S_1$, if neither $x$ nor $y$ have an edge to $S_2$, then we get at least $\lceil \alpha n \rceil$ induced copies of $H$ in $G$ such that $V(H) \cap S_1 = \{x, y\}$ and $V(H) \cap S_2 \neq \emptyset$. Likewise, the same is true for each edge in $S_2$. Let $h(G)$ denote the number of induced copies of $H$ in $G$.

Each edge with one end in $S_1$ and one end in $S_2$ will be incident to at most $|S_1|+|S_2|$ edges contained in either $S_1$ or $S_2$. Therefore, each edge $xy$ of the at most $2 t_1 n - |S_1|-|S_2|$ edges between $S_1$ and $S_2$ is part of at most $|S_1|+|S_2|$ triples $\{x, y, z\}$ such that $z \in S_1 \cup S_2$ but $\{x, y, z\}$ does not induce a copy of $H$. It follows that $h(G)$ is at least:
\begin{align}
    &2 \lceil \alpha n \rceil \left((1-\epsilon)\frac{t_1^2 n^2}{2} - \frac{(t_1 n - \lceil \alpha n \rceil)^2}{2} -\frac{t_1 n - \lceil \alpha n \rceil}{2}\right) 
    \\ &- \left(2t_1 n - |S_1| - |S_2|\right)\left(|S_1| + |S_2|\right). \label{eq:q}
\end{align}
Using the fact that $|S_1| \ge \alpha n > \frac{t_1 n}{2}$ and likewise for $|S_2|$, we conclude that replacing each of $|S_1|, |S_2|$ by $\lceil \alpha n \rceil$ does not increase \eqref{eq:q}, and so $h(G)$ is at least: 
\begin{align*}
        &2 \lceil \alpha n \rceil \left((1-\epsilon)\frac{t_1^2 n^2}{2} - \frac{(t_1 n - \lceil \alpha n \rceil)^2}{2}\right) \\
        &- \lceil \alpha n \rceil (t_1 n -\lceil \alpha n \rceil) - 2 \lceil \alpha n \rceil \left(2t_1 n - 2 \lceil \alpha n \rceil \right),
\end{align*}
which equals:
\begin{align}
    2 \lceil \alpha n \rceil \left((1-\epsilon)\frac{t_1^2 n^2}{2} - \frac{(t_1 n - \lceil \alpha n \rceil)^2}{2}\right) - 5 \lceil \alpha n \rceil (t_1 n -\lceil \alpha n \rceil). \label{eq:alphanew}
\end{align}

Observing again that $\alpha > t_1 / 2$, we obtain that:
\begin{align*}
    h(G) \ge 2 \alpha n \left((1-\epsilon)\frac{t_1^2 n^2}{2} - \frac{(t_1 n - \alpha n)^2}{2}\right) - 5( t_1 n  - \alpha n)(\alpha n).
\end{align*}
Here we used that $(1-\epsilon)\frac{t_1^2 n^2}{2} - \frac{(t_1 n - \alpha n)^2}{2} > 0$ since $\epsilon < 1/2$ and $(1-\alpha/t_1) < 1/2$. Therefore, changing the first occurrence $\lceil \alpha n \rceil $ to $\alpha n$ in \eqref{eq:alphanew} does not increase its value. For the second occurrence, notice that $0 \leq t_1n- \lceil \alpha n \rceil \leq t_1n - \alpha n$, and so $- \frac{(t_1 n - \alpha n)^2}{2} \leq - \frac{(t_1 n - \lceil \alpha n \rceil)^2}{2}$. For the last two occurrences, we notice that $\alpha n > t_1n/2$, and $- 5( t_1 n  - x)x$ is a convex function in $x$ minimized when $x= t_1n/2$.

Let us simplify:
\begin{align*}
    h(G) &\ge 2 \alpha n \left((1-\epsilon)\frac{t_1^2 n^2}{2} - \frac{(t_1 n - \alpha n)^2}{2}\right) - 5( t_1 n  - \alpha n)\alpha n\\
    &= \alpha n \left((t_1n)^2 (1 - \epsilon - (1-\alpha/t_1)^2) - 5(t_1-\alpha)n\right)
\end{align*}

We apply Theorem \ref{triange_count}, using that $m = (1+\delta)t_1 n^2$. This tells us that the number of triangles in $G$ is at least:
\begin{align*}
& \frac{h(G)}{3} + \frac{4m}{3n}\left(m-\frac{n^2}{4}\right) \\
\geq \ & 
    \frac{\alpha n}{3} \left((t_1n)^2 (1 - \epsilon - (1-\alpha/t_1)^2) - 5(t_1-\alpha)n\right)\\
    + & \frac{4(1+\delta_1)t_1n^3 }{3}\left((1+\delta_1)t_1-\frac{1}{4}\right).
\end{align*}

Since there is no rainbow triangle, then each triangle contains at least two edges of the same color, and it follows that the number of triangles is at most:
\begin{align*}
    (1+\delta_1)n \binom{t_1 n}{2} = (1+\delta_1)n\frac{t_1 n(t_1 n -1)}{2}.
\end{align*}
We claim that the following is true, which will give an immediate contradiction:
\begin{align}
\begin{split}    &\frac{\alpha n}{3} \left((t_1n)^2 (1 - \epsilon - (1-\alpha/t_1)^2) - 5(t_1-\alpha)n\right) \\
    &+ \frac{4(1+\delta_1)t_1n^3}{3}\left((1+\delta_1)t_1-\frac{1}{4}\right)\\
>& \ (1+\delta_1)n\frac{t_1 n(t_1 n-1)}{2}. 
\end{split} \label{eq:thebigone}
\end{align}
%Note that we may assume that $5t_1 \epsilon < 1$ because $\epsilon \le 1/2$ and $t_1 < 2/5$ since otherwise we are done by Theorem \ref{existing}. Then the assumption that $5t_1 \epsilon < 1$ implies that:
%\begin{align*}
%    \frac{5}{3}\alpha(t_1-\alpha) &< 5\alpha(t_1-\alpha) \\
%    &= \frac{5t_1^2}{4}(1+\sqrt{1-2\epsilon})(1-\sqrt{1-2\epsilon}) \\
%    &= \frac{5t_1^2 \epsilon}{2} \\
%    &< \frac{t_1}{2} \\
%    &< \frac{t_1(1+\delta_1)}{2},
%\end{align*}
Note that 
$$\frac{5}{3}\alpha(t_1-\alpha) \leq \frac{5t_1}{12} < \frac{t_1(1+\delta_1)}{2},$$
and so \eqref{eq:thebigone} is implied by the following (after cancelling terms corresponding to the above inequality and multiplying by $n^3/6$): 
\begin{align*}
\begin{split}    & 2\alpha  \left(t_1^2 (1 - \epsilon - (1-\alpha/t_1)^2) \right) \\
    &+ (1+\delta_1)t_1\left(8(1+\delta_1)t_1-2\right)\\
>& \ 3(1+\delta_1)t_1^2. 
\end{split} 
\end{align*}
Now, let's plug in $\alpha = \frac{t_1}{2}(1+\lambda)$ where $\lambda = \sqrt{1-2\epsilon}$; and then divide by $t_1$; then we would like to show that \eqref{eq:themediumone} holds.
\begin{align}
\begin{split}    & (1 + \lambda)  \left(t_1^2 (1 - \epsilon - (1-(1+\lambda)/2)^2) \right) \\
    &+ (1+\delta_1)\left(8(1+\delta_1)t_1-2\right)\\
&- 3(1+\delta_1)t_1 > 0. 
\end{split} \label{eq:themediumone}
\end{align}

We claim that the derivative of the above expression with respect to $t_1$ and with respect to $\delta_1$ is positive for all $(t_1,\delta_1)$ with $t_1 \ge t$ and $\delta_1 \ge \delta$. Indeed, the derivative of the left-hand side of \eqref{eq:themediumone} with respect to $t_1$ is equal to:
\begin{align*}
     &(1+\lambda)2t_1(1-\varepsilon-(1-\lambda)^2/4) + (1+\delta_1)(8(1+\delta_1)-3)\\
    > & (1+\delta_1)(8(1+\delta_1)-3) > 0
\end{align*}
since $\lambda < 1$ and $\epsilon < 1/2$ and $\delta_1 > 0$. 

The derivative of the left-hand side of \eqref{eq:themediumone} with respect to $\delta_1$ is equal to:
\begin{align*}
    16t_1(1+\delta_1) - 2 - 3t_1 > 0,
\end{align*}
since by \eqref{eq:213} we have:
\begin{align*}
    \frac{16}{3}(1+\delta_1) \ge \frac{16}{3}(1+\delta) > \frac{2}{3t}+1 \ge \frac{2}{3t_1}+1.
\end{align*}

Recall that by \eqref{eq:212}, 
    $$\frac{8}{3t}\left(t(1+\delta)-\frac{1}{4}\right)(1+\delta)+\frac{t}{12}(\lambda+1)\left(4(1-\epsilon)-\left(1-\lambda\right)^2\right) - (1+\delta) > 0.
    $$
Multiplying by $3t$ gives: 
$$\left(8t(1+\delta)-2\right)(1+\delta)+\frac{t^2}{4}(\lambda+1)\left(4(1-\epsilon)-\left(1-\lambda\right)^2\right) - 3t(1+\delta) > 0,$$
which is equal to \eqref{eq:themediumone} when $t = t_1$ and $\delta = \delta_1$. Since the left-hand-size is increasing in both of these variables, it follows that \eqref{eq:themediumone} holds. This gives the desired contradiction, and completes the proof. \end{proof}

This immediately gives Theorem \ref{main1}, which we restate here:
%\begin{thm}
%$(1.1077,1/3)$ is triangular; $(1.3481,1/4)$ is triangular.
%\end{thm}
\mainone*
\begin{proof}
One can verify that $t=1/3$, $\delta=0.1077$, and $\epsilon=0.4746$ satisfy the above inequalities, and also that $t=1/4$, $\delta = 0.3481$, and $\epsilon = 0.2774$ also satisfy the above inequalities. This completes the proof.\end{proof}

In general, for a fixed value of $t$ we determine the minimum value of $\delta$ that this result gives by using the following tool, which we link here: \\\url{https://www.desmos.com/calculator/1x06qadpqp}.

\section{$n$ colors of size at least $t n$}
In this section, we show Theorem \ref{main2} by first showing the following more general statement.
\begin{thm}\label{thm_3}
For real $0 < \epsilon < 1/2$, $0 < \delta \le 1$, $t > 0$, suppose that the following conditions hold:
\begin{align}
     \frac{(1-\epsilon)t-\delta}{1-\delta} \ge 0.3465; \label{eq:chbound} \\
t > \frac{1}{2}-(1-\epsilon)^2 t^2; \label{eq:t}\\
 \frac{8}{3}\left(t-\frac{1}{4}\right)>(1-\delta(2\epsilon - 2\epsilon^2))t; \label{eq:83} \\
1+2\delta \epsilon^2 > 4\epsilon \delta. \label{eq:delta}
\end{align}
If $G$ is a simple edge-colored graph with $n$ colors classes, each of size at least $tn$, then there exists a rainbow triangle in $G$.
\end{thm}
\begin{proof}
Suppose that $G$ does not contain a rainbow triangle. Let $t_1$ be a real number such that $t_1 n = \lceil tn \rceil$. We may assume that each color class has size exactly equal to $t_1 n$. Now, let us call a color $c$ \emph{concentrated} if there exists a vertex which is incident to at least $(1-\epsilon)t_1n$ edges of color $c$. 

    \begin{claim}\label{claim:delta}
        More than $(1-\delta)n$ colors are concentrated. 
    \end{claim}

Suppose not; that is, at least $\delta n$ colors are not concentrated. For each such color $c$, the number of triangles with at least two edges of color $c$ is at most the number of happy triples in a graph with $t_1n$ edges and maximum degree at most $\lfloor (1-\epsilon)t_1n \rfloor$. Let us set $l = \lfloor (1-\epsilon)t_1n \rfloor$ unless $\lfloor (1-\epsilon)t_1n \rfloor < t_1n/2$; set $l = (t_1n+1)/2$ in this case instead. Note that since $\epsilon < 1/2$, the second case only happens if $t_1n$ is odd and $\lfloor (1-\epsilon)t_1n \rfloor = (t_1n-1)/2$. Therefore, for each color which is not concentrated, its edges form a graph of maximum degree at most $l$. The number of triangles in $G$ is at most the sum of the number of happy triples in each color, because $G$ has no rainbow triangle, and so each triangle is a happy triple for some color. 

Therefore, using Lemma \ref{lem:ht}, the number of triangles in $G$ is at most 
\begin{align}
    (1-\delta)n {t_1n \choose 2} + \delta n \left({l \choose 2} + {t_1n-l+1 \choose 2} \right). \label{eq:19}
\end{align}
Using Lemma \ref{lem:convex}, we know that $f_{t_1n}(l) \leq f_{t_1n}((1-\varepsilon)t_1n)$ unless $l < (t_1n+1)/2$, because $(t_1n+1)/2$ is the global minimum of $f_{t_1n}$. This leaves only the case that $l = t_1n/2$, which we consider separately. 

\paragraph{Case 1: $l \neq t_1n/2$} 
Then, by \eqref{eq:19} and using $f_{t_1n}(l) \leq f_{t_1n}((1-\varepsilon)t_1n)$, the number of triangles in $G$ is at most 
\begin{align*}
 &(1-\delta)n {t_1n \choose 2} + \delta n\left({(1-\varepsilon)t_1n \choose 2} + {t_1n - (1-\varepsilon)t_1n + 1 \choose 2}\right)\\
 = &(1-\delta)n {t_1n \choose 2} + \frac{\delta n}{2}\left(((1-\varepsilon)t_1n)^2 + (t_1n - (1-\varepsilon)t_1n)^2 + t_1n(2 \varepsilon - 1)\right).
 \end{align*}
We would like to add $(1-\delta)n^2t_1/2 - \delta n^2t_1\frac{2 \epsilon - 1}{2}$ to this to cancel terms not divisible by $t_1^2n^3$. So we claim that $(1-\delta)n^2t_1/2 - \delta n^2t_1\frac{2 \epsilon - 1}{2} \geq 0$. It suffices to show that $(1-\delta) - \delta (2 \epsilon - 1) \geq 0$. This holds since $\epsilon < 1/2$, so $2\epsilon - 1 < 0$; and $\delta < 1$.

Therefore, the number of triangles in $G$ is at most 
\begin{align*}
 & (1-\delta)n \frac{(t_1n)^2}{2} + \frac{\delta n}{2}\left(((1-\varepsilon)t_1n)^2 + (t_1n - (1-\varepsilon)t_1n)^2\right)\\
=\ &\frac{t_1^2n^3}{2}((1-\delta) + \delta((1-\varepsilon)^2 + \epsilon^2))\\
=\ &\frac{t_1^2n^3}{2}(1 - \delta(2\epsilon - 2 \varepsilon^2)). 
 \end{align*}
 By Theorem \ref{goodman_thm}, we have that the number of triangles is at least:
\begin{align*}
    \frac{4m}{3n}\left(m-\frac{n^2}{4}\right) \ge \frac{4t_1 n}{3}\left(t_1 n^2-\frac{n^2}{4}\right).
\end{align*}
Comparing the two bounds for the number of triangles in $G$, we have shown that 
\begin{align}
    \frac{4}{3}\left(t_1 - \frac{1}{4}\right) \leq \frac{t_1}{2}(1 - \delta(2\epsilon - 2 \varepsilon^2)). \label{eq:t1}
\end{align}
Comparing the derivatives of both sides with respect to $t_1$, we notice that decreasing $t_1$ preserves the inequality (since $(1 + \delta(2\epsilon^2 - 2 \varepsilon)) \leq 1$). But $t \leq t_1$, and so \eqref{eq:t1} holds for $t$ instead of $t_1$, contrary to \eqref{eq:83}. This finishes Case 1. 
 
\paragraph{Case 2: $l = t_1n/2$} 
Then, by \eqref{eq:19}, the number of triangles in $G$ is at most 
\begin{align*}
 &(1-\delta)n {t_1n \choose 2} + \delta n (t_1n/2)^2\\
 =& (1-\delta)n ((t_1n)^2/2 - t_1n/2)+ \delta n (t_1n/2)^2
 \end{align*}
 and again at least 
\begin{align*}
    \frac{4m}{3n}\left(m-\frac{n^2}{4}\right) \ge \frac{4t_1 n}{3}\left(t_1 n^2-\frac{n^2}{4}\right).
\end{align*}
by Theorem \ref{goodman_thm}. Comparing the two bounds for the number of triangles, we obtain that 
$$(1-\delta)n ((t_1n)^2/2 - t_1n/2)+ \delta n (t_1n/2)^2 \geq \frac{4t_1 n}{3}\left(t_1 n^2-\frac{n^2}{4}\right).$$
Let us simplify this by cancelling $t_1n^2/4$: 
\begin{align}
    2(1-\delta) (t_1n - 1)+ \delta t_1n = 2\delta - 2 - \delta t_1n + 2t_1n\geq \frac{4}{3}\left(4t_1 n-n\right). \label{eq:whew}
\end{align}

We would like to show that 
\begin{align}
    2\delta - 2 - \delta t_1n + 2t_1n \leq 2t_1n(1 - \delta(2\epsilon - 2 \varepsilon^2)).\label{eq:yay}
\end{align} 
Cancelling gives: 
$$ 2\delta - 2 - \delta t_1n\leq - 2t_1n(\delta(2\epsilon - 2 \varepsilon^2)).$$
Therefore, we want to prove:  
$$ 2\delta - 2 \leq t_1n\delta(1 - 4\epsilon + 4 \varepsilon^2) = \delta t_1n(1-2\epsilon)^2.$$
The left-hand side is negative, and the right-hand side is non-negative, and so the desired inequality \ref{eq:yay} holds. Together with \eqref{eq:whew}, we conclude that
 $$2t_1n(1 - \delta(2\epsilon - 2 \varepsilon^2)) \geq \frac{4}{3}\left(4t_1 n-n\right).$$
Dividing by $4n$, we have once again shown that \eqref{eq:t1} holds, which gives a contradiction as before. This concludes case 2, and therefore the proof of Claim \ref{claim:delta}.

\medskip

It follows that at least $(1-\delta)n$ of the colors are concentrated. Now, there are two cases. Suppose first that for each vertex $v \in V(G)$, there is at most one color concentrated at $v$. Here, we proceed with a method inspired by the reduction of Conjecture \ref{ch_thm} to Conjecture \ref{aharoni_thm} found in \cite{aharoni}. Let $T \subseteq V(G)$ be the set of vertices $v$ for which there exists a color $c_v$ concentrated at $v$. Form a digraph $D$ with vertex set equal to $T$, and for each vertex $v$ and edge $e=vu$ of color $c_v$ which is incident to $v$, add an arc from $v$ to $u$. It follows that $D$ has minimum out-degree at least $(1-\epsilon)t_1 n-\left(n-|T|\right)$, and therefore since $|T| \ge (1-\delta)n$, we have that:
\begin{align*}
    \frac{\delta^+(D)}{|D|} \ge \frac{(1-\epsilon)t_1 n-\left(n-|T|\right)}{|T|} \ge \frac{(1-\epsilon)t_1-\delta}{1-\delta}. 
\end{align*}

By \eqref{eq:chbound}, we have that:
$$\frac{(1-\epsilon)t_1-\delta}{1-\delta} \ge \frac{(1-\epsilon)t-\delta}{1-\delta} \ge 0.3465,$$
so it follows from Theorem \ref{ch_result_thm} that there exists a directed triangle $T$ in $D$. Looking at the underlying edges in $G$ corresponding to the directed edges of $T$, we obtain a rainbow triangle in $G$.

\medskip

So, we may assume that there is some vertex $v$ at which two colors $c_1$ and $c_2$ are concentrated. For $i \in \{1, 2\}$, let $S_i$ be the set of vertices $u$ such that $vu$ is an edge of color $c_i$. Then, since there are no rainbow triangles, it follows that the all edges with one end in $S_1$ and one end in $S_2$ have color either $c_1$ or $c_2$. Therefore, we get at least $(1-\epsilon)^2 t_1^2 n^2 - 2\epsilon t_1 n$ non-edges in $G$, and thus the number of edges is at most:
\begin{align*}
    \frac{n(n-1)}{2}-(1-\epsilon)^2 t_1^2 n^2 +  2\epsilon t_1 n \le \frac{n^2}{2}-(1-\epsilon)^2 t_1^2 n^2,
\end{align*}
since $2\epsilon t_1 < 1/2$ because $\epsilon < 1/2$ and $t_1 < 0.4$ (otherwise we are done by Theorem \ref{existing}). However, we know that there are at least $t_1 n^2$ edges in $G$, and therefore, 
\begin{align}\label{eq:finally}
    t_1 \leq \frac{1}{2} - (1-\epsilon)^2t_1^2
\end{align} 
Decreasing $t_1$ decreasing the right-hand side and increases the left-hand side; so \eqref{eq:finally} also holds for $t$ instead of $t_1$, contrary to \eqref{eq:t}. This is a contradiction. Therefore, $G$ contains a rainbow triangle, as desired. This completes the proof.\end{proof}

Now, we immediately obtain Theorem \ref{main2}, which we restate here:
%\begin{thm}
%$(1,0.3988)$ is triangular.
%\end{thm}
\maintwo*
\begin{proof}
We verify that $\epsilon = 0.03846$, $\delta = 0.0681$, and $t = 0.3988$ satisfy all the conditions of Theorem \ref{thm_3}. Since the edges of $G$ are colored with $n$ colors each of size at least $tn=0.3988n$, it follows that there exists a rainbow triangle, as desired. This completes the proof.\end{proof}

\section{Concluding Remarks and Further Work}
In this paper, we found improved bounds for the triangle case of Conjecture \ref{aharoni_thm}. We found our methods to be significantly more effective in the case of having at least $(1+\delta)n$ colors, which matches our intuition that this case is significantly easier. Further work could improve on either of our results, but we find the prospect of a new method introducing significant progress to the case of $n$ colors to be particularly interesting.

\section{Acknowledgments}

We acknowledge the support of the Natural Sciences and Engineering Research Council of Canada (NSERC), [funding reference number RGPIN-2020-03912]. Cette recherche a été financée par le Conseil de recherches en sciences naturelles et en génie du Canada (CRSNG), [numéro de référence RGPIN-2020-03912]. This research was funded in part by the Government of Ontario. 

We are thankful to Michelle Delcourt and Brittany Pittman for pointing out an error in Lemma \ref{triange_count} in an earlier version  of this paper. We thank the referees for their helpful comments.

\appendix

\section{Computation for $r = 3$ and $r = 4$} \label{ap:1}

    Below we show the calculation for plugging $r = 4$ into Equation \ref{eq:k}. 
    \begin{align*}
        2\left(\frac{k}{4} + 6 - \frac{7}{2}\right)\left(\frac{k}{4} + 6 - \frac{9}{2}\right) + 3k/2 &> k^2/4 - 1/4\\
        \left(k + 24 - 14\right)\left(k + 24 - 18\right) + 12k &> 2k^2 - 2\\
        k^2 + 28k + 60 &> 2k^2 - 2\\
        0 &> k^2 - 28k - 62\\
        14 + \sqrt{196 + 62} &> k\\
        30.07 &> k.
    \end{align*}

    And for $r = 3$: 
    \begin{align*}
        \frac{3}{2}\left(\frac{k}{3} + \frac{9}{2} - \frac{7}{2}\right)\left(\frac{k}{3} + \frac{9}{2} - \frac{9}{2}\right) + 12k &> k^2/4 - 1/4\\
        \left(2k + 27 - 21\right)k + 12k &> 3k^2 - 3\\
        2k^2 + 18k &> 3k^2 - 3\\
        0 &> k^2 - 18k - 3\\
        9 + \sqrt{81+3} &> k\\
        18.17 &> k.
    \end{align*}

\section{The case $k \leq 103$} \label{ap:2}

Our goal is to prove Lemma \ref{lem:ht}, which we restate:  
\begin{lemma}Let $k, l \in \mathbb{N}$. Let $G$ be a graph with $k$ edges and with maximum degree at most $l \geq k/2$. Then $G$ contains at most ${l \choose 2} + {k-l+1 \choose 2}$ happy triples. 
\end{lemma} 
We already proved: 
\begin{itemize}
    \item It suffices to consider the case $k \leq 103$.
    \item We may assume every vertex has degree less than $k/2$. 
    \item As a function $l$, the expression ${l \choose 2} + {k-l+1 \choose 2}$ is increasing for $l \geq (k+1)/2$ and equal at $k/2$ and $k/2+1$. 
\end{itemize}
Therefore, it suffices to consider the case $l = \lceil k/2 \rceil$: every graph $G$ for which we have not proved the lemma yet has maximum degree less than $k/2$, and so is valid for every $l \geq k/2$, and the upper bound we aim to prove is minimized at $l = \lceil k/2 \rceil$ among all $l \in \{\lceil k/2 \rceil,\lceil k/2 \rceil + 1, \dots, k\}$. 

\begin{table}[]
    \centering
    \begin{tabular}{c|c||c|c||c|c||c|c||c|c}
         $k, l$ & bound & $k, l$ & bound & $k, l$ & bound & $k, l$ & bound & $k, l$ & bound \\
         \hline
         $3, 2$ & $2$ & $4, 2$ & $4$ & $5, 3$ & $6$ & $6, 3$ & $9$ & $7, 4$ & $12$\\ 
$8, 4$ & $16$ & $9, 5$ & $20$ & $10, 5$ & $25$ & $11, 6$ & $30$ & $12, 6$ & $36$\\ 
$13, 7$ & $42$ & $14, 7$ & $49$ & $15, 8$ & $56$ & $16, 8$ & $64$ & $17, 9$ & $72$\\ 
$18, 9$ & $81$ & $19, 10$ & $90$ & $20, 10$ & $100$ & $21, 11$ & $110$ & $22, 11$ & $121$\\ 
$23, 12$ & $132$ & $24, 12$ & $144$ & $25, 13$ & $156$ & $26, 13$ & $169$ & $27, 14$ & $182$\\ 
$28, 14$ & $196$ & $29, 15$ & $210$ & $30, 15$ & $225$ & $31, 16$ & $240$ & $32, 16$ & $256$\\ 
$33, 17$ & $272$ & $34, 17$ & $289$ & $35, 18$ & $306$ & $36, 18$ & $324$ & $37, 19$ & $342$\\ 
$38, 19$ & $361$ & $39, 20$ & $380$ & $40, 20$ & $400$ & $41, 21$ & $420$ & $42, 21$ & $441$\\ 
$43, 22$ & $462$ & $44, 22$ & $484$ & $45, 23$ & $506$ & $46, 23$ & $529$ & $47, 24$ & $552$\\ 
$48, 24$ & $576$ & $49, 25$ & $600$ & $50, 25$ & $625$ & $51, 26$ & $650$ & $52, 26$ & $676$\\ 
$53, 27$ & $702$ & $54, 27$ & $729$ & $55, 28$ & $756$ & $56, 28$ & $784$ & $57, 29$ & $812$\\ 
$58, 29$ & $841$ & $59, 30$ & $870$ & $60, 30$ & $900$ & $61, 31$ & $930$ & $62, 31$ & $961$\\ 
$63, 32$ & $992$ & $64, 32$ & $1024$ & $65, 33$ & $1056$ & $66, 33$ & $1089$ & $67, 34$ & $1122$\\ 
$68, 34$ & $1156$ & $69, 35$ & $1190$ & $70, 35$ & $1225$ & $71, 36$ & $1260$ & $72, 36$ & $1296$\\ 
$73, 37$ & $1332$ & $74, 37$ & $1369$ & $75, 38$ & $1406$ & $76, 38$ & $1444$ & $77, 39$ & $1482$\\ 
$78, 39$ & $1521$ & $79, 40$ & $1560$ & $80, 40$ & $1600$ & $81, 41$ & $1640$ & $82, 41$ & $1681$\\ 
$83, 42$ & $1722$ & $84, 42$ & $1764$ & $85, 43$ & $1806$ & $86, 43$ & $1849$ & $87, 44$ & $1892$\\ 
$88, 44$ & $1936$ & $89, 45$ & $1980$ & $90, 45$ & $2025$ & $91, 46$ & $2070$ & $92, 46$ & $2116$\\ 
$93, 47$ & $2162$ & $94, 47$ & $2209$ & $95, 48$ & $2256$ & $96, 48$ & $2304$ & $97, 49$ & $2352$\\ 
$98, 49$ & $2401$ & $99, 50$ & $2450$ & $100, 50$ & $2500$ & $101, 51$ & $2550$ & $102, 51$ & $2601$\\ 
$103, 52$ & $2652$ \\
    \end{tabular}
    \caption{Outputs for $k \in \{3, \dots, 109\}$ and $l = \lceil k/2 \rceil$.}
    \label{tab:my_label}
\end{table}

We proceed via dynamic programming, and include the code in this appendix (written in Python). The output of the code for the cases we are interested in, that is, $k \in \{3, \dots, 103\}$ and $l = \lceil k/2 \rceil$, is shown in Table \ref{tab:my_label}. We have checked by computer that the output matches the upper bound that we aim to show in each case. 

Let us describe how the program works. For fixed $k$ and $l$ (dropping the condition that $l \geq k/2$), we would like to find the maximum number of happy triples in a graph $G$ with $k$ edges and maximum degree at most $l$. When $l = 1$, we have no happy triples. When $l = 2$, we have that each edge $xy$ is part of at most two happy triples (because each of $x, y$ has at most one more incident edge). Since each happy triple involves at least two edges, it follows that when $l=2$, the number of happy triples for $k$ edges is at most $k$. When $k = 3$, the maximum is actually 2, because $G = K_3$ has only one happy triple; so a 3-edge path actually achieves the maximum of 2. We create an array that stores an upper bound for maximum number of happy triples in a graph with $k$ edges and maximum degree at most $l$ for each pair $(k, l)$. We initialize the answer for $(k, l)$ with $l = 2$ and $k \in \{4, \dots, 103\}$ as $k$. 

Now, let us fix a pair $(k, l)$ with $k, l \in \mathbb{N}$ with $2 \leq k, l \leq 103$. Let us assume that we know the array values for all pairs $(k', l')$ with $k' < k$ and $l' \leq l$. 

As shown in the proof of Lemma \ref{lem:ht}, when $l \geq k$, this maximum is achieved by $K_{1,k}$ with a value of ${k \choose 2}$. So we may assume that $l < k$. Let $G$ be a graph achieving the maximum number of happy triples for the pair $(k, l)$. Then, we consider a vertex $v$ of maximum degree $j$ in $G$; so $j \in \{1, \dots, l\}$. Let $F$ be the set of edges of $G$ incident with $v$, and let $F'$ be the set of remaining edges. Then the number of happy triples in $G$ is at most the sum of: 
\begin{itemize}
    \item the number of happy triples containing two edges in $F$, which is ${j \choose 2}$; 
    \item the number of happy triples containing $v$ and an edge in $F'$, which is at most $k-j = |F'|$; 
    \item the number of happy triples not containing $v$, and therefore only using edges in $F'$, with the additional property that $(V(G), F')$ has maximum degree at most $j$ (from the choice of $v$), which is at most the array entry at $(k-j, j)$. 
\end{itemize}
Therefore, we compute this number for each $j$, and take the maximum over all choices of $j$; this is assigned as the array value for $(k, l)$. The above shows that the array value for $(k, l)$ is an upper bound for the number of happy triples in a graph with $k$ edges and maximum degree at most $l$. 

\begin{python}
# bound we want to prove
def bound(k, l):                                                         
    return l*(l-1)/2 + (k-l+1)*(k-l)/2
    
# range for k: we consider k between 1 and maximum-1
maximum = 200
    
# initialize and all-0 array for the dynamic program
# ar[k][l] will upper bound the number of happy triples when we 
# have k edges and maximum degree at most l
ar = [[0]*maximum for i in range(maximum)]

# initialize some values in the array that we know
for k in range(4, maximum): 
    ar[k][2] = k 

# in the dynamic program, to compute ar[k][l], we need to know 
# a[k'][l'] for certain values that have l' <= l and k' < k
# this order of processing entries accomplishes this
for l in range(2, maximum):
    for k in range(2, maximum):
        # do not overwrite our preset values
        if ar[k][l] > 0: continue
        # if l >= k, then the best solution is K_{1, k}, which 
        # gives {k choose 2} happy triples
        if k <= l: 
            ar[k][l] = k*(k-1)/2
            continue
        # otherwise, consider a vertex v of maximum degree in an 
        # optimum solution; it has degree 1 <= j <= l
        bid = 0
        goodj = -1
        for j in range(1, l+1): 
            # then we get {j choose 2} happy triples from edges 
            # incident with j
            # each edge not among those j forms at most 1 happy 
            # triple with v, for at most k-j overall
            # there are k-j remaining edges, and their maximum 
            # degree is at most j
            happytriples = j*(j-1)/2 + k-j + ar[k-j][j]
            bid = max(bid, happytriples)
            if bid == happytriples: goodj = j
        # assign best choi/ce over all values of j
        ar[k][l] = bid
        # if we outperform the bound, output this (adding +0.2 for
        # possible floating point issues)
        if bid >= bound(k, l)+0.2: 
            print(k,l,bid,goodj, bound(k, l), ar[k-j][j])
\end{python}

\end{document}